\documentclass[a4paper,12pt]{amsart}
\usepackage{amsfonts}
\usepackage{amssymb}
\usepackage{ifthen}
\usepackage{graphicx}

\usepackage{color}
\newtheorem{thm}{Theorem}
\newtheorem{cor}[equation]{Corollary}
\newtheorem{lem}{Lemma}
\newtheorem{prop}[equation]{Proposition}

\newtheorem{conj}[equation]{Conjecture}
\newtheorem{rem}{Remark}
\theoremstyle{definition}
\newtheorem{defn}{Definition}%[section]
\newtheorem{examp}{Example}
\newtheorem{prob}[equation]{Problem}
\newtheorem{ques}[equation]{Question}

\newcounter {own}
\def\theown {\thesection       .\arabic{own}}

\newenvironment{pf}[1][]{%
 \vskip 1mm
 \noindent
 \ifthenelse{\equal{#1}{}}%
  {{\slshape Proof. }}%
  {{\slshape #1.} }%
 }%
{\qed\smallskip}

\newcounter{alphabet}
\newcounter{tmp}

\makeatletter
\newcommand{\Ref}[1]{\@ifundefined{r@#1}{}{\setcounter{tmp}{\ref{#1}}\Alph{tmp}}}
\makeatother
\newenvironment{Lem}[1][]{\refstepcounter{alphabet}%
\bigskip%
\noindent%
{\bf Lemma \Alph{alphabet}}%
{\bf .} \itshape}{\vskip 8pt}

\newcommand{\IR}{{\mathbb R}}

\newcommand{\ID}{{\mathbb D}}
\newcommand{\IB}{{\mathbb B}}

\newcommand{\Aut}{{\operatorname{Aut}}}

\def\be{\begin{equation}}
\def\ee{\end{equation}}

\newcommand{\bee}{\begin{enumerate}}
\newcommand{\eee}{\end{enumerate}}

\newcommand{\blem}{\begin{lem}}
\newcommand{\elem}{\end{lem}}
\newcommand{\bthm}{\begin{thm}}
\newcommand{\ethm}{\end{thm}}
\newcommand{\bcor}{\begin{cor}}
\newcommand{\ecor}{\end{cor}}
\newcommand{\beg}{\begin{examp}}
\newcommand{\eeg}{\end{examp}}
\newcommand{\begs}{\begin{examples}}
\newcommand{\eegs}{\end{examples}}
\newcommand{\bdefe}{\begin{defn}}
\newcommand{\edefe}{\end{defn}}
\newcommand{\bprob}{\begin{prob}}
\newcommand{\eprob}{\end{prob}}
\newcommand{\bques}{\begin{ques}}
\newcommand{\eques}{\end{ques}}
\newcommand{\bei}{\begin{itemize}}
\newcommand{\eei}{\end{itemize}}

\newcommand{\bcon}{\begin{conj}}
\newcommand{\econ}{\end{conj}}
\newcommand{\bcons}{\begin{conjs}}
\newcommand{\econs}{\end{conjs}}
\newcommand{\bprop}{\begin{prop}}
\newcommand{\eprop}{\end{prop}}
\newcommand{\br}{\begin{rem}}
\newcommand{\er}{\end{rem}}
\newcommand{\brs}{\begin{rems}}
\newcommand{\ers}{\end{rems}}
\newcommand{\bo}{\begin{obser}}
\newcommand{\eo}{\end{obser}}
\newcommand{\bos}{\begin{obsers}}
\newcommand{\eos}{\end{obsers}}
\newcommand{\bpf}{\begin{pf}}
\newcommand{\epf}{\end{pf}}
\newcommand{\ba}{\begin{array}}
\newcommand{\ea}{\end{array}}
\newcommand{\beq}{\begin{eqnarray}}
\newcommand{\beqq}{\begin{eqnarray*}}
\newcommand{\eeq}{\end{eqnarray}}
\newcommand{\eeqq}{\end{eqnarray*}}

\newcommand{\ds}{\displaystyle}

\newcounter{minutes}
\divide\time by 60
\newcounter{hours}
\multiply\time by 60 \addtocounter{minutes}{-\time}
\begin{document}
%\begin{center}
%{ Communicated by  }
%\end{center}
\bibliographystyle{amsplain}
\title [] {Weighted Lipschitz continuity, Schwarz-Pick's Lemma and Landau-Bloch's theorem
for hyperbolic-harmonic mappings in  $\mathbb{C}^{n}$ }

%%%%%%%% BEGIN TIMESTAMP
\def\thefootnote{}
\footnotetext{ \texttt{\tiny File:~\jobname .tex,
          printed: \number\day-\number\month-\number\year,
          \thehours.\ifnum\theminutes<10{0}\fi\theminutes}
} \makeatletter\def\thefootnote{\@arabic\c@footnote}\makeatother
%%%%%%%% END TIMESTAMP

\author{SH. Chen}
\address{Sh. Chen, Department of Mathematics,
Hunan Normal University, Changsha, Hunan 410081, People's Republic
of China.} \email{shlchen1982@yahoo.com.cn}

\author{S. Ponnusamy}
\address{S. Ponnusamy, Department of Mathematics,
Indian Institute of Technology Madras, Chennai-600 036, India.}
\email{samy@iitm.ac.in}
\author{X. Wang${}^{~\mathbf{*}}$}
\address{X. Wang, Department of Mathematics,
Hunan Normal University, Changsha, Hunan 410081, People's Republic
of China.} \email{xtwang@hunnu.edu.cn}

 \keywords{Hyperbolic-harmonic mapping,  Bloch
spaces, Landau-Bloch's theorem,
Schwarz-Pick's Lemma.\\
${}^{\mathbf{*}}$ Corresponding author}
%\date{\today  %November 4, 10; File: Ch-W-S13-${}_{}$hyper-harm${}_{}$submit.tex}
\begin{abstract}
In this paper, we discuss some properties on hyperbolic-harmonic
mappings in the unit ball of $\mathbb{C}^{n}$. First, we investigate
the relationship between the weighted Lipschitz functions and the
hyperbolic-harmonic Bloch spaces. Then we establish the Schwarz-Pick
type theorem for hyperbolic-harmonic mappings and apply it to prove
the existence of Landau-Bloch constant for  mappings in
 $\alpha$-Bloch spaces.
\end{abstract}

\thanks{The research was partly supported by NSF of China (No. 11071063). }

\maketitle \pagestyle{myheadings} \markboth{SH. Chen,  S. Ponnusamy,
X. Wang}{On hyperbolic-harmonic mappings}

\section{Introduction and preliminaries }\label{csw-sec1}

Let $\mathbb{C}^{n}$ denote the complex Euclidean $n$-space.
For $z=(z_{1},\ldots,z_{n})\in \mathbb{C}^{n}$,
the conjugate of $z$, denoted by $\overline{z}$, is defined by
$\overline{z}=(\overline{z}_{1},\ldots, \overline{z}_{n} ). $ For
$z$ and $w=(w_{1},\ldots,w_{n})\in\mathbb{C}^{n}
$, the standard Hermitian scalar product on $\mathbb{C}^{n}$ and the Euclidean
norm of $z$ are given by
$$\langle z,w\rangle := \sum_{k=1}^nz_k\overline{w}_k\quad \mbox{and}\quad
|z|:={\langle z,z\rangle}^{1/2}=(|z_{1}|^{2}+\cdots+|z_{n}|^{2})^{ 1/2},
$$
respectively. For $a\in \mathbb{C}^n$,  $\IB^n(a,r)=\{z\in \mathbb{C}^{n}:\, |z-a|<r\}$ is the (open) ball of radius $r$ with center $a$.
Also, we let $\IB^n(r):=\IB^n(0,r)$ and use $\IB^n$ to denote the unit ball $\IB^n(1)$, and
$\mathbb{D}=\mathbb{B}^1$.
%, and for $a\in\mathbb{R}^{n}$, $\IB_{R}^n(a, r)=\{x\in \mathbb{R}^{n}:\, |x-a|<r\}.$
We can interpret $ \mathbb{C}^n$ as the real $2n$-space $ \mathbb{R}^{2n}$ so that
a ball in $ \mathbb{C}^n$ is also a ball in $ \mathbb{R}^{2n}$.
Following the standard convention, for $a\in\mathbb{R}^{n}$, we may let
$\IB_{\IR}^n(a, r)=\{x\in \mathbb{R}^{n}:\, |x-a|<r\}$ so that
$\IB_{\IR}^n(r):=\IB_{\IR}^n(0,r)$ and $\mathbb{B}^{n}_{\mathbb{R}}=\IB_{\IR}^n(1)$
denotes the open unit ball in $\mathbb{R}^{n}$ centered at the origin.

\begin{defn}\label{defn0}
A twice continuously differentiable complex-valued function
$f=u+iv$ on $\mathbb{B}^{n}$ is called a
{\it hyperbolic-harmonic} (briefly, $h$-harmonic, in the following) if and only if the real-valued
functions $u$ and $v$ satisfy $\Delta_{h}u=\Delta_{h}v=0$ on $\mathbb{B}^{n}$, where
$$\Delta_{h}:=(1-|z|^{2})^{2}\sum_{k=1}^{n}\left(\frac{\partial}{\partial x^{2}_{k}}
+\frac{\partial}{\partial
y^{2}_{k}}\right)+4(n-1)(1-|z|^{2})\sum_{k=1}^{n}\left(x_{k}\frac{\partial}{\partial
x_{k}} +y_{k}\frac{\partial}{\partial y_{k}}\right)
$$
denotes the {\it Laplace-Beltrami operator} and $z_{k}=x_{k}+iy_{k}$ for $k=1,\ldots,n$.
\end{defn}

Obviously, when
$n=1$, all $h$-harmonic mappings are planar harmonic
mappings. We refer to \cite{B1,EO,GJ,TW} for more details of
$h$-harmonic mappings.

%The {\it hyperbolic Poisson kernel} applied to solve the Dirichlet
%problem for $\Delta_{h}$ is defined by
%$$\mbox{P}_{h}(z,\zeta)=\left(\frac{1-|z|^{2}}{|z-\zeta|^{2}}\right)^{2n-1}
%$$
%for $z\in\mathbb{B}^{n}$ and $\zeta\in\partial\mathbb{B}^{n}$.
It turns out that if $\psi\in C(\partial\mathbb{B}^{n})$, then the Dirichlet problem
$$\begin{cases}
\displaystyle \Delta_{h}f=0
& \mbox{ in } \mathbb{B}^{n}\\
\displaystyle f=\psi &\mbox{ on }\, \partial\mathbb{B}^{n}
\end{cases}
$$
has unique solution in $C(\overline{\mathbb{B}}^{n})$ and can be represented by
$$f(z)=\int_{\partial\mathbb{B}^{n}}\mbox{P}_{h}(z,\zeta)\psi(\zeta)d\sigma(\zeta),
$$
where $d\sigma$ is the unique normalized surface measure on
$\partial\mathbb{B}^{n}$ and $\mbox{P}_{h}(z,\zeta)$ is
the {\it hyperbolic Poisson kernel} defined by
$$\mbox{P}_{h}(z,\zeta)=\left(\frac{1-|z|^{2}}{|z-\zeta|^{2}}\right)^{2n-1}
\quad \mbox{($z\in\mathbb{B}^{n}$, $\zeta\in\partial\mathbb{B}^{n}$).}
$$
Here $C(\Omega)$ stands for the set of all continuous functions on $\Omega$.
A planar harmonic mapping $f$ in $\mathbb{D}$ is called a {\it
harmonic Bloch mapping} if and only if
$$\beta_{f}=\sup_{ z,w\in\mathbb{D},\ z\neq w}\frac{|f(z)-f(w)|}{\rho(z,w)}<\infty.
$$
Here $\beta_{f}$ is the {\it  Lipschitz number} of $f$ and
$$\rho(z,w)=\frac{1}{2}\log\left(\frac{1+|\frac{z-w}{1-\overline{z}w}|}
{1-|\frac{z-w}{1-\overline{z}w}|}\right)=\mbox{arctanh}\Big|\frac{z-w}{1-\overline{z}w}\Big|
$$
denotes the hyperbolic distance between $z$ and $w$ in $\mathbb{D}$. It can be proved that
$$\beta_{f}=\sup_{z\in\mathbb{D}}\big\{(1-|z|^{2})\big[|f_{z}(z)|+|f_{\overline{z}}(z)|\big]\big\}.
%\ (\mbox{see also \cite{ CPW2,CPW3}}).
$$
We refer to \cite[Theroem 2]{Co} (see also \cite{CPW2,CPW3}) for a proof of the last fact.

For a complex-valued $h$-harmonic mapping $f$ on $\mathbb{B}^{n}$, we introduce
%the following notations:
$$\widehat{\nabla f}=\left(\frac{\partial
f}{\partial z_{1}},\ldots,\frac{\partial f}{\partial z_{n}}\right
)\;\;\mbox{and}\;\; \widehat{\nabla\overline{f} }=\left
(\frac{\partial f}{\partial \overline{z}_{1}},\ldots,\frac{\partial
f}{\partial \overline{z}_{n}}\right ).
$$

\begin{defn}\label{defn1}
The {\it $h$-harmonic Bloch space $\mathcal{HB}$} consists
of complex-valued $h$-harmonic mappings $f$ defined on $\mathbb{B}^{n}$ such that
$$\|f\|_{\mathcal{HB}}=\sup_{z\in\mathbb{B}^{n}}\big\{(1-|z|^{2})\big[|\widehat{\nabla f}(z) |+|
\widehat{\nabla \overline{f}}(z) |\big]\big\}<\infty.
$$
\end{defn}
Obviously, when  $n=1$, $\|f\|_{\mathcal{HB}}=\beta_f$.
\noindent For a pair of distinct points $z$ and $w$ in $\mathbb{B}^{n}$, let
$$\mathcal{L}_{f}(z,w)=\frac{(1-|z|^{2})^{\frac{1}{2}}\big(1-|w|^{2}\big)^{\frac{1}{2}}|f(z)-f(w)|}{|z-w|}
$$
denote the {\it weighted Lipschitz function} of a given $h$-harmonic
mapping $f:\,\mathbb{B}^{n}\rightarrow \mathbb{C}$. The relationship
between weighted Lipschitz functions and (analytic) Bloch spaces has
attracted much attention (cf. \cite{MVM,Co, HW,L,MVM1,Pav}). Our
first aim is to characterize the mappings in $h$-harmonic Bloch
spaces in terms of their corresponding weighted Lipschitz functions.
This is done in Theorem \ref{thm1s} which is indeed a generalization
of \cite[Theorem 1]{Co} and \cite[Theorem 3]{HW}.

Throughout, $\mathcal{H}(\IB^n, \mathbb{C}^n)$ denotes the set of all continuously differentiable
mappings $f$ from $\mathbb{B}^{n}$ into $\mathbb{C}^{n}$ with $f=(f_{1},\ldots,f_{n})$ and
$f_{i}(z)=u_{j}(z)+iv_{j}(z)$, where $u_{j}$ and $v_{j}$  are real-valued mappings on $\mathbb{B}^{n}$.
For $f\in \mathcal{H}(\IB^n, \mathbb{C}^n)$,
%  with $f_{j}(z)=u_{j}(z)+iv_{j}(z)$ for $j=1,\ldots,n$.
%,  $u_{j}$ and $v_{j}$  are real-valued mappings on $\mathbb{B}^{n}$. Here $z=(z_{1},\ldots,z_{n})$ with
%$z_{k}=x_{k}+iy_{k}$, $k\in\{1,\ldots,n\}$,  so that
the real Jacobian matrix of $f$ is given by
$$J_{f}=\left(\begin{array}{cccc}
\ds \frac{\partial u_{1}}{\partial x_{1}}\;~~ \frac{\partial
u_{1}}{\partial y_{1}}\;~~ \frac{\partial u_{1}}{\partial x_{2}}\;~~
\frac{\partial u_{1}}{\partial y_{2}}\;~~\cdots\;~~
\frac{\partial u_{1}}{\partial x_{n}}\;~~ \frac{\partial u_{1}}{\partial y_{n}}\\[2mm]
\ds \frac{\partial v_{1}}{\partial x_{1}}\;~~ \frac{\partial v_{1}}{\partial y_{1}}\;~~
\frac{\partial v_{1}}{\partial x_{2}}\;~~ \frac{\partial
v_{1}}{\partial y_{2}}\;~~\cdots\;~~
 \frac{\partial v_{1}}{\partial x_{n}}\;~ \frac{\partial v_{1}}{\partial
 y_{n}}\\[2mm]
\ds \vdots \hspace{1cm} \vdots \\[2mm]
 \ds \frac{\partial u_{n}}{\partial x_{1}}\;~ \frac{\partial u_{n}}{\partial y_{1}}\;~~
\frac{\partial u_{n}}{\partial x_{2}}\; \frac{\partial
u_{n}}{\partial y_{2}}\;~~\cdots\;~
 \frac{\partial u_{n}}{\partial x_{n}}\;~~ \frac{\partial u_{n}}{\partial y_{n}}\\[2mm]
\ds \frac{\partial v_{n}}{\partial x_{1}}\;~~ \frac{\partial v_{n}}{\partial y_{1}}\;~~
\frac{\partial v_{n}}{\partial x_{2}}\;~~ \frac{\partial
v_{n}}{\partial y_{2}}\;~~\cdots\;~~
 \frac{\partial v_{n}}{\partial x_{n}}\;~~ \frac{\partial v_{n}}{\partial
 y_{n}}
\end{array}\right).
$$

A vector-valued mapping $f\in\mathcal{H}(\IB^n, \mathbb{C}^n)$ is said to be {\it $h$-harmonic},
if each component $f_{i}$ ($1\leq i\leq n$) is $h$-harmonic
mapping from $\mathbb{B}^{n}$ into $\mathbb{C}$. We denote by $\mathcal{H}_{h}(\mathbb{B}^{n},
\mathbb{C}^n)$ the set of all vector-valued $h$-harmonic mappings from $\mathbb{B}^{n}$ into
$\mathbb{C}^n$.

For each $f=(f_{1},\ldots,f_{n})\in \mathcal{H}(\IB^n,
\mathbb{C}^n)$, denote by
$$f_{z}=\big(\widehat{\nabla f}_1,\ldots,\widehat{\nabla f}_n\big)^T
$$
the matrix formed by the complex gradients $\widehat{\nabla
f}_1,\ldots,\widehat{\nabla f}_n$, and let
$$f_{\overline{z}}=\big(\widehat{\nabla\overline{f}}_1,\ldots,\widehat{\nabla\overline{f}}_n\big)^T,
$$
where $T$ means the matrix transpose.

For an $n\times n$ matrix $A=\left (a_{ij}\right )_{n\times n}$, the
operator norm of $A$ is given by
$$|A|=\sup_{z\neq 0}\frac{|Az|}{|z|}=\max\big\{|A\theta|:\,
\theta\in\partial \mathbb{B}^{n}\big\}.
$$
Then for $f \in \mathcal{H}(\IB^n, \mathbb{C}^n)$, we use the standard
notations:
\be\label{eqbe1}
\Lambda_{f}(z)=\max_{
\theta\in\partial\mathbb{B}^{n}}|f_{z}(z)\theta+f_{\overline{z}}(z)\overline{\theta}|\;\;
\mbox{ and }\;\; \lambda_{f}(z)=\min_{
\theta\in\partial\mathbb{B}^{n}}|f_{z}(z)\theta+f_{\overline{z}}(z)\overline{\theta}|.
\ee

We see that (see for instance \cite{HG1})
\be\label{eqbe2}
\Lambda_{f}=\max_{\theta\in\partial\mathbb{B}^{2n}_{\mathbb{R}}}|J_{f}\theta|\;\;
\mbox{ and }\;\;
\lambda_{f}=\min_{\theta\in\partial\mathbb{B}^{2n}_{\mathbb{R}}}|J_{f}\theta|.
\ee

Let $\mathcal{PH}(\IB^n, \mathbb{C}^n)$ denote the set of all
$f=(f_{1},\ldots,f_{n})\in \mathcal{H}(\IB^n, \mathbb{C}^n)$ such
that all partial derivatives $\partial f_{j}/z_{k}$ and
$\partial f_{j}/\overline{z}_{k}$ $(1\leq j,k\leq n)$ are $h$-harmonic in $\IB^n$.

We remark that when $n=1$, every harmonic mapping from $\mathbb{D}$ to $\mathbb{C}$
belongs to $\mathcal{PH}(\mathbb{D}, \mathbb{C})$. The converse is not true as the
function $f(z)=|z|^2$ shows.

\begin{defn}\label{defn4}
For $\alpha>0,$ the {\it vector-valued  $h$-harmonic $\alpha$-Bloch space
$\mathcal{HB}_{n}(\alpha)$} consists of all mappings in
$\mathcal{PH}(\IB^n, \mathbb{C}^n)$ such that
$$\|f\|_{\mathcal{HB}_n(\alpha)}=\sup_{z\in\mathbb{B}^{n}}\big\{(1-|z|^{2})^{\alpha}\big[| f_{z}(z) |+|
f_{\overline{z}}(z) |\big]\big\}<\infty.
$$
\end{defn}

Obviously, $\mathcal{HB}_{1}(\alpha)$ contains the
harmonic $\alpha$-Bloch space as a proper subset (see \cite{CPW3}).
One of the long standing open problems in function theory is to
determine the precise value of the univalent Landau-Bloch constant
for analytic functions of $\ID$. In recent years, this problem has attracted much
attention, see \cite{Bo,LM, M1} and references therein. For general holomorphic
mappings of more than one complex variable, no Landau-Bloch constant
exists (cf. \cite{W}). In order to obtain some analogs of
Landau-Bloch's theorem for mappings with several complex variables,
it became necessary to restrict the class of mappings considered (cf.
\cite{B,HG1,CPW4,LX,Ra,T,W}).

Based on Heinz's Lemma and Colonna's distortion theorem (\cite[Theorem 3]{Co}) for
planar harmonic mappings, in \cite{HG1}, the authors established
the Schwarz-Pick type theorem  for bounded pluriharmonic mappings and pluriharmonic
$K$-mappings. As a consequence of it, the authors in \cite{HG1} obtained
Landau-Bloch theorem as generalizations of the main results from
\cite{HG}. It is known that every pluriharmonic mapping $f$ defined in $\mathbb{B}^{n}$
admits a decomposition $f=h+\overline{g}$, where $h$ and $g$ are holomorphic in
$\mathbb{B}^{n}$. This decomposition property is no longer valid for mappings
in $\mathcal{HB}_{n}(\alpha)$. Hence the methods of proof used in \cite{HG1} are no
longer applicable for mappings in $\mathcal{H}_{h}(\mathbb{B}^{n}, \mathbb{C}^n)$ and $\mathcal{HB}_{n}(\alpha)$.
In view of this reasoning, in this article, we use entirely a different approach and
prove Schwarz-Pick type theorem for mappings in $\mathcal{H}_{h}(\mathbb{B}^{n},
\mathbb{C}^n)$ and then establish the Landau-Bloch theorem for
mappings in $\mathcal{HB}_{n}(\alpha)$ (see Theorems \ref{thm1} and \ref{thm3w}).
It is worth pointing out that Theorems \ref{thm1} and \ref{thm3w} are indeed
generalizations of \cite[Theorem 1]{Co} and \cite[Theorem 2.4]{CPW3}, respectively. % \cite[Theorem 2]{HG}.

\section{Characterization of mappings in $h$-harmonic Bloch spaces }\label{csw-sec2}

Consider the group $\Aut(\mathbb{B}^{n})$ consisting of all biholomorphic
mappings of $\mathbb{B}^{n}$ onto itself. Then for each $a\in\mathbb{B}^{n}$,
$\phi_{a}$ defined by \cite{R}:
$$\phi_{a}(z)=\frac{a-P_{a}z-(1-|a|^{2})^{\frac{1}{2}} (z-P_{a}z)}{1-\langle z,a\rangle}
$$
belongs to $\Aut(\mathbb{B}^{n})$, where
%$\langle z,a\rangle=z_{1}\overline{a}_{1}+\cdots+z_{n}\overline{a}_{n},$
$P_{a}z=\frac{a\langle z, a\rangle}{\langle a, a\rangle}$. Moreover, we find that
\be\label{eq1}
1-|\phi_{a}(z)|^{2}=\frac{(1-|z|^{2})(1-|a|^{2})}{|1-\langle
z,a\rangle|^{2}}.
\ee

Using the arguments as in the proof of \cite[Lemma 2.5]{MV}, we have

\begin{lem}\label{lem1f}
Suppose $f:\,\overline{\mathbb{B}}_{\IR}^{n}(a,r)\rightarrow\mathbb{R}$ is a
continuous, and $h$-harmonic in $\mathbb{B}_{\IR}^{n}(a,r)$. Then
$$|\nabla f(a)|\leq
\frac{2(n-1)\sqrt{n}}{nV(n)r^{n}}\int_{\partial\mathbb{B}_{\IR}^{n}(a,r)}|f(a)-f(t)|\,d\sigma(t),
$$
where $d\sigma$ denotes the surface measure on
$\partial\mathbb{B}_{\IR}^{n}(a,r)$ and $V(n)$, the volume of the unit
ball in $\mathbb{R}^{n}$.
\end{lem}
\bpf Without loss of generality, we may assume that $a=0$ and
$f(0)=0.$ Let
$$K(x,t)=\frac{1}{nr^{n-1}V(n)}\left(\frac{r^{2}-|x|^{2}}{|x-t|^{2}}\right)^{n-1}.
$$
Then by the assumption on $f$, we see that \cite{B1}
$$f(x)=\int_{\partial\mathbb{B}_{\IR}^{n}(r)}K(x,t)f(t)\,d\sigma(t),
~\mbox{ $x\in\mathbb{B}_{\IR}^{n}(r)$}.
$$
Further, a
%, where  $d\sigma$ denotes the surface measure on $\partial\mathbb{B}_{\IR}^{n}(a,r)$.
computation shows that
$$\frac{\partial}{\partial x_{i}}K(x,t)=\frac{-2(n-1)\big(r^{2}-|x|^{2}\big)^{n-2}}{nr^{n-1}V(n)}\cdot
\frac{\big[|x-t|^{2}x_{i}+(r^{2}-|x|^{2})(x_{i}-t_{i})\big]}{|x-t|^{2n}}
$$
which yields
$$\frac{\partial}{\partial x_{i}}K(0,t)=\frac{2(n-1)t_{i}}{nV(n)r^{n+1}}
$$
whence
\begin{eqnarray*}
|\nabla f(0)|&=&\left[\sum_{i=1}^{n}\left|\int_{\partial\mathbb{B}_{\IR}^{n}(r)}\frac{\partial}{\partial
x_{i}}K(0,t)f(t)\,d\sigma(t)\right |^{2}\right]^{\frac{1}{2}}\\
&\leq&\sum_{i=1}^{n}\left|\int_{\partial\mathbb{B}_{\IR}^{n}(r)}\frac{\partial}{\partial
x_{i}}K(0,t)f(t)\,d\sigma(t)\right |\\
&\leq&\int_{\partial\mathbb{B}_{\IR}^{n}(r)}|f(t)|\sum_{i=1}^{n}\Big|\frac{\partial}{\partial
x_{i}}K(0,t)\Big|\,d\sigma(t)\\
&\leq&\sqrt{n}\int_{\partial\mathbb{B}_{\IR}^{n}(r)}|f(t)|\left (\sum_{i=1}^{n}\Big|\frac{\partial}{\partial
x_{i}}K(0,t)\Big|^{2}\right )^{\frac{1}{2}}\,d\sigma(t)\\
&=&\frac{2(n-1)\sqrt{n}}{nV(n)r^{n}}\int_{\partial\mathbb{B}_{\IR}^{n}(r)}|f(t)|\,d\sigma(t),
\end{eqnarray*}
from which the lemma follows.
\epf

\begin{lem}\label{lem2x}
Let $f=u+iv$ be a continuously differentiable mapping from
$\mathbb{B}^{n}$ into $\mathbb{C}$, where $u$ and $v$ are
real-valued functions. Then for $z\in\mathbb{B}^{n}$,
\be\label{eqs1} |\widehat{\nabla f}(z)|+|\widehat{\nabla
\overline{f}}(z)|\leq |\nabla u(z)|+|\nabla v(z)|, \ee where $\nabla
u=\Big( \frac{\partial u}{\partial x_{1}}, \frac{\partial
u}{\partial y_{1}},\cdots, \frac{\partial u}{\partial x_{n}},
\frac{\partial u}{\partial y_{n}}\Big)$ and $\nabla v=\Big(
\frac{\partial v}{\partial x_{1}}, \frac{\partial v}{\partial
y_{1}},\cdots, \frac{\partial v}{\partial x_{n}}, \frac{\partial
v}{\partial y_{n}} \Big)$.
\end{lem} \bpf
By the basic change of variables, for each $k=1,2,\ldots,n,$ we have
$$f_{z_{k}}(z)=\frac{1}{2}\big(f_{x_{k}}(z)-if_{y_{k}}(z)\big)~
\mbox{and}~f_{\overline{z}_{k}}(z)=\frac{1}{2}\big(f_{x_{k}}(z)+if_{y_{k}}(z)\big)
$$
which implies
$$f_{z_{k}}(z)=\frac{1}{2}\big[u_{x_{k}}(z)+v_{y_{k}}(z)+i(v_{x_{k}}(z)-u_{y_{k}}(z))\big]
$$
and similarly,
$$f_{\overline{z}_{k}}(z)=\frac{1}{2}\big[u_{x_{k}}(z)-v_{y_{k}}(z)+i(v_{x_{k}}(z)+u_{y_{k}}(z))\big].
$$
The classical Cauchy-Schwarz inequality gives
$$ |\widehat{\nabla
f}(z)|=\frac{1}{2}\sqrt{\sum_{k=1}^{n}\Big[\big(u_{x_{k}}(z)+v_{y_{k}}(z)\big)^{2}+\big(v_{x_{k}}(z)-u_{y_{k}}(z)\big)^{2}\Big]}
\leq\frac{1}{2}\big(|\nabla u(z)|+|\nabla v(z)|\big)
$$
and similarly,
$$|\widehat{\nabla \overline{f}}(z)|=\frac{1}{2}\sqrt{\sum_{k=1}^{n}\Big[\big(u_{x_{k}}(z)-v_{y_{k}}(z)\big)^{2}
+\big(v_{x_{k}}(z)+u_{y_{k}}(z)\big)^{2}\Big]}
\leq\frac{1}{2}\big(|\nabla u(z)|+|\nabla v(z)|\big)
$$
from which we obtain the desired inequality \eqref{eqs1}.
\epf

%We remark that the equality sign in (\ref{eqs1}) does not always hold. This can be seen from the following example.

\begin{examp}
Consider $f(z)=z^{2}+\overline{z}=u(x,y)+iv(x,y)$ so that $u(x,y)=x^{2}+x-y^{2}$ and $v(x,y)=2xy-y$.
It is easy to see that
$$|f_{z}(0)|+|f_{\overline{z}}(0)|=1 ~\mbox{ and }~ |\nabla u(0)|+|\nabla v(0)|=2
$$
showing that strict inequality in (\ref{eqs1}) is possible.
\end{examp}

%\begin{thm}\label{thm1s}
%Let $f$ be an $h$-harmonic mapping in $\mathbb{B}^n$.
%$f\in \mathcal{HB}$ if and only if
%$\ds \sup_{z,w\in\mathbb{B}^{n},z\neq w}\mathcal{L}_{f}(z,w)<\infty.
%$
%\Big\{(1-|z|^{2})^{\frac{1}{2}}(1-|w|^{2})^{\frac{1}{2}}\frac{|f(z)-f(w)|}{|z-w|}\Big\}$$
%where $E(z,r)=\{w\in\mathbb{B}^{n}:\, |\phi_{z}(w)|<r\ \mbox{and}\
%\phi_{z}(w)\in\Aut(\mathbb{B}^{n})\}.$
%\end{thm}

\begin{thm}\label{thm1s}
%Let $f$ be an $h$-harmonic mapping in $\mathbb{B}^n$.
$f\in \mathcal{HB}$ if and only if $\ds
\sup_{z,w\in\mathbb{B}^{n},z\neq w}\mathcal{L}_{f}(z,w)<\infty. $
%\Big\{(1-|z|^{2})^{\frac{1}{2}}(1-|w|^{2})^{\frac{1}{2}}\frac{|f(z)-f(w)|}{|z-w|}\Big\}$$
%where $E(z,r)=\{w\in\mathbb{B}^{n}:\, |\phi_{z}(w)|<r\ \mbox{and}\
%\phi_{z}(w)\in\Aut(\mathbb{B}^{n})\}.$
\end{thm}
\bpf  First we prove the  necessity. For each pair of distinct
points $z$ and $w$ in $\mathbb{B}^{n}$, we have
\begin{eqnarray*}
|f(z)-f(w)|&=&\left |\int_{0}^{1}\frac{df}{dt}(zt+(1-t)w)\, dt\right |\\
&=&\left |\sum_{k=1}^{n}(z_{k}-w_{k})\int_{0}^{1}\frac{df}{d\varsigma_{k}}(zt+(1-t)w)\, dt \right .\\
&&\hspace{1cm}\left .+\sum_{k=1}^{n}(\overline{z}_{k}-\overline{w}_{k})\int_{0}^{1}\frac{df}{d\overline{\varsigma}_{k}}(zt+(1-t)w)\,dt\right|\\
&\leq&\sum_{k=1}^{n}|z_{k}-w_{k}|\cdot\left |\int_{0}^{1}\frac{df}{d\varsigma_{k}}(zt+(1-t)w)\, dt\right|\\
&&\hspace{1cm}+\sum_{k=1}^{n}|\overline{z}_{k}-\overline{w}_{k}|\cdot\left
|\int_{0}^{1}\frac{df}{d\overline{\varsigma}_{k}}(zt+(1-t)w)
\,dt\right |,\\
\end{eqnarray*} where
$\varsigma=(\varsigma_{1},\cdots,\varsigma_{n})=zt+(1-t)w$. Hence we
see that
\begin{eqnarray*}
|f(z)-f(w)|&\leq&\left(\sum_{k=1}^{n}|z_{k}-w_{k}|^{2}\right)^{\frac{1}{2}}
\left\{\left [\sum_{k=1}^{n}\left (\int_{0}^{1}\Big|\frac{\partial
f}{\partial \varsigma_{k}}(zt+(1-t)w)\Big|\,dt\right )^{2}\right]^{\frac{1}{2}}\right .\\
&&\hspace{2cm} \left . +\left [\sum_{k=1}^{n}\left
(\int_{0}^{1}\Big|\frac{\partial
f}{\partial \overline{\varsigma}_{k}}(zt+(1-t)w)\Big|\, dt\right )^{2}\right ]^{\frac{1}{2}}\right \} \\
%&\leq& \sqrt{n}|z-w|\left [\int_{0}^{1} |\nabla f(tz+(1-t)w)|\,dt
%\right . \\ &&
%\hspace{2cm}  \left . +\int_{0}^{1}|\nabla \overline{f}(tz+(1-t)w)|\,dt\right ].
&\leq& \sqrt{n}|z-w|\int_{0}^{1} \left [ |\widehat{\nabla
f}(tz+(1-t)w)| + |\widehat{\nabla \overline{f}}(tz+(1-t)w)|\right
]dt.
\end{eqnarray*}
%It follows from letting $\psi(t)=tz+(1-t)w$ that
This gives
\begin{eqnarray*}
\frac{|f(z)-f(w)|}{|z-w|}&\leq&\sqrt{n}\int_{0}^{1}\frac{\big[|\widehat{\nabla
f}(\varsigma)|+|\widehat{\nabla
\overline{f}}(\varsigma)|\big](1-|\varsigma|^{2})}{1-|\varsigma|^{2}}\,dt\\
&\leq&\sqrt{n}\|f\|_{\mathcal{HB}}\int_{0}^{1}\frac{dt}{1-|\varsigma|^{2}}\\
&\leq&\sqrt{n}\|f\|_{\mathcal{HB}}\int_{0}^{1}\frac{dt}
{\big[(1-t)(1-|z|)\big]^{\frac{1}{2}}\big[t(1-|w|)\big]^{\frac{1}{2}}}\\
&=&\frac{\pi\sqrt{n}\|f\|_{\mathcal{HB}}}{(1-|z|)^{\frac{1}{2}}(1-|w|)^{\frac{1}{2}}}.
\end{eqnarray*}
%where $\psi(t)=tz+(1-t)w$.
Thus,
$$\sup_{z,w\in\mathbb{B}^{n},z\neq w}\mathcal{L}_{f}(z,w)\leq \pi\sqrt{n}\|f\|_{\mathcal{HB}}.
$$

Next we prove the sufficiency part. Let $f=u+iv,$ where $u$ and $v$
are real $h$-harmonic functions. Fix $r\in (0,1)$. In view of
(\ref{eq1}) and the fact that $|\langle z,a\rangle| \leq |z|\,|a|$,
we easily have \be\label{eq2t}
|\phi_{a}(z)|\leq\frac{|z-a|}{|1-\langle
z,a\rangle|}\leq\frac{|z-a|}{1-|a|}, \ee
whence for $a\in\mathbb{B}^{n}$, %by (\ref{eq2}), we have
$$\mathbb{B}^{n}\Big(a,\frac{r(1-|a|^{2})}{2}\Big)\subset E(a,r),
$$
where $E(a,r)=\{z\in\mathbb{B}^{n}:\, |\phi_{a}(z)|<r\}.$ By Lemma
\ref{lem1f}, we have
\begin{eqnarray*}
(1-|z|^{2})|\nabla u(z)|&\leq&\frac{(2n-1)\sqrt{2n}(1-|z|^{2})}{n
V(2n)[\frac{r(1-|z|^{2})}{2}]^{2n}}
\int_{\partial\mathbb{B}^{n}\big(z,\frac{r(1-|z|^{2})}{2}\big)}|u(\zeta)-u(z)|\,d\sigma(\zeta)\\
&=&M(|z|,r)\int_{\partial\mathbb{B}^{n}
\big(z,\frac{r(1-|z|^{2})}{2}\big)}|u(\zeta)-u(z)|\,d\sigma(\zeta),
\end{eqnarray*}
where $V(2n)$ denotes the volume of the unit ball in
$\mathbb{R}^{2n}~(\mbox{or}~\mathbb{C}^{n})$ and
$$M(|z|,r)=\frac{2^{2n}(2n-1)\sqrt{2n}}{nV(2n)(1-|z|^{2})^{2n-1}r^{2n}}.
$$ Similarly, we obtain
$$(1-|z|^{2})|\nabla v(z)|\leq M(|z|,r)\int_{\partial\mathbb{B}^{n}
\big(z,\frac{r(1-|z|^{2})}{2}\big)}|v(\zeta)-v(z)|\,d\sigma(\zeta).
$$
%Cauchy's inequality and chain rules of derivation show that
%$$|\nabla f(z)|\leq\frac{1}{2}(|\nabla u(z)|+|\nabla v(z)|)~\mbox{
%and }~|\nabla \overline{f}(z)|\leq\frac{1}{2}(|\nabla u(z)|+|\nabla
%v(z)|),$$ which implies
By Lemma \ref{lem2x}, we have
\begin{eqnarray*}
(1-|z|^{2})(|\widehat{\nabla f}(z)|+|\widehat{\nabla
\overline{f}}(z)|)&\leq&(1-|z|^{2})(|\nabla u(z)|+|\nabla v(z)|)\\
&\leq&M(|z|,r)\int_{\partial\mathbb{B}^{n}\big(z,\frac{r(1-|z|^{2})}{2}\big)}
\Big(|u(\zeta)-u(z)|\\
&& \hspace{2cm}  +|v(\zeta)-v(z)|\Big)\,d\sigma(\zeta)\\
&\leq&\sqrt{2}M(|z|,r)M_{1}\int_{\partial\mathbb{B}^{n}\big(z,\frac{r(1-|z|^{2})}{2}\big)}\,d\sigma(\zeta)\\
&=&\frac{4\sqrt{n}(2n-1)}{r}M_{1},
\end{eqnarray*}
where $M_{1}=\sup\{|f(z)-f(w)|:\, w\in E(z,r)\}.$ Hence for all
$w\in\mathbb{B}^{n}\Big(z,\frac{r(1-|z|^{2})}{2}\Big)\subset
E(z,r),$ it follows from (\ref{eq1}) and (\ref{eq2t}) that
\begin{eqnarray*}
\frac{(1-|z|^{2})^{\frac{1}{2}}(1-|w|^{2})^{\frac{1}{2}}}{|z-w|}
&=&\frac{(1-|z|^{2})^{\frac{1}{2}}
(1-|w|^{2})^{\frac{1}{2}}}{|1-\langle z,w\rangle|}\cdot\frac{|1-\langle z,w\rangle|}{|z-w|}\\
&=&\sqrt{1-|\phi_{z}(w)|^{2}}\cdot\frac{|1-\langle z,w\rangle|}{|z-w|}\\
&\geq&\sqrt{1-r^{2}}\cdot\frac{|1-\langle z,w\rangle|}{|z-w|}\\
&\geq&\frac{\sqrt{1-r^{2}}}{r}.
\end{eqnarray*}
Therefore, there exists a positive constant $M_{2}(n,r)$ such that
\begin{eqnarray*}
(1-|z|^{2})[|\widehat{\nabla f}(z) |+|\widehat{\nabla
\overline{f}}(z) |] \leq M_{2}(n,r)\sup_{w\in E(z,r),w\neq
z}\mathcal{L}_{f}(z,w),
\end{eqnarray*}
from which we see that $f\in\mathcal{HB}.$
%The proof is complete.
\epf

\section{ Schwarz-Pick type theorem and Landau-Bloch theorem}\label{csw-sec3}

 The following result is a Schwarz-Pick type theorem
for $h$-harmonic mappings in $\mathcal{H}_{h}(\mathbb{B}^{n}, \mathbb{C}^n)$.

\begin{thm}\label{thm1}
Let $f\in \mathcal{H}_{h}(\mathbb{B}^{n}, \mathbb{C}^n)$ with $|f(z)|\leq M$ for
$z\in\mathbb{B}^n$, where $M$ is a positive constant. Then
\be\label{eq-1}
\left|f(z)-\frac{(1-|z|)^{2n-1}}{(1+|z|)^{2n-1}}f(0)\right|\leq
M\left[1-\frac{(1-|z|)^{2n-1}}{(1+|z|)^{2n-1}}\right]
\ee
and
\be\label{eq2}
\Lambda_{f}\leq \frac{2(2n-1)M}{(1-|z|)^{2}}.
\ee
\end{thm}
\bpf We first prove (\ref{eq-1}). Without loss of generality, we assume that
$f$ is also $h$-harmonic on $\partial\mathbb{B}^{n}.$ The hyperbolic Poisson integral
formula states that
\be\label{eq-x1}
f(z)=\int_{\partial\mathbb{B}^{n}}\mbox{P}_{h}(z,\zeta)f(\zeta)\,d\sigma(\zeta),
\quad \int_{\partial\mathbb{B}^{n}}\mbox{P}_{h}(z,\zeta)\,d\sigma(\zeta) =1.
\ee%=
%\int_{\partial\mathbb{B}^{n}}\frac{(1-|z|^{2})^{2n-1}}{|z-\zeta|^{2(2n-1)}}f(\zeta)\,d\sigma(\zeta).$$
%where $d\sigma$ denotes the normalized surface measure on
%$\partial\mathbb{B}^{n}$.
%where $$ $$
As $\mbox{P}_{h}(0,\zeta)=1$ and $\mbox{P}_{h}(z,\zeta)|\leq 1$ for $\zeta \in \partial\mathbb{B}^{n}$ and all
$z\in \mathbb{B}^{n}$, the representation \eqref{eq-x1} quickly yields that
\begin{eqnarray*}
\left|f(z)-\frac{(1-|z|)^{2n-1}}{(1+|z|)^{2n-1}}f(0)\right|&=&\left|\int_{\partial\mathbb{B}^{n}}
\left [\frac{(1-|z|^{2})^{2n-1}}{|z-\zeta|^{2(2n-1)}}-\frac{(1-|z|)^{2n-1}}{(1+|z|)^{2n-1}}\right]f(\zeta)
\,d\sigma(\zeta)\right|\\
&\leq&\int_{\partial\mathbb{B}^{n}}
\left[\frac{(1-|z|^{2})^{2n-1}}{|z-\zeta|^{2(2n-1)}}-\frac{(1-|z|)^{2n-1}}{(1+|z|)^{2n-1}}\right ]|f(\zeta)|
\,d\sigma(\zeta)\\
&\leq&M\left[1-\frac{(1-|z|)^{2n-1}}{(1+|z|)^{2n-1}}\right]
\end{eqnarray*}
and the proof of (\ref{eq-1}) follows.

Next, we  prove (\ref{eq2}). Let $f=(f_{1},\ldots,f_{n})$ and
$\theta=(\theta_{1},\ldots,\theta_{n})^{T}
\in\partial\mathbb{B}^{n}$. Without loss of generality, we assume
that $f$ is also $h$-harmonic on $\partial\mathbb{B}^{n}.$
%By
%the Poisson integral formula, we have
%$$f(z)=\int_{\partial\mathbb{B}^{n}}\frac{(1-|z|^{2})^{2n-1}}{|z-\zeta|^{2(2n-1)}}f(\zeta)\,d\sigma(\zeta),
%$$
%where $d\sigma$ denotes the normalized surface measure on$\partial\mathbb{B}^{n}$.
If we consider the formula (\ref{eq-x1}) for $f$ componentwise and then the partial derivatives
with respect to the variables $z_k$ and $\overline{z}_{k}$, we see that
$$(f_{j}(z))_{z_{k}}=\int_{\partial\mathbb{B}^{n}}
\frac{-(2n-1)(1-|z|^{2})^{2n-2}[\overline{z}_{k}|\zeta-z|^{2}+(1-|z|^{2})(\overline{z}_{k}-\overline{\zeta}_{k})]}
{|z-\zeta|^{4n}}f_{j}(\zeta)\,d\sigma(\zeta)
$$
and
$$(f_{j}(z))_{\overline{z}_{k}}=\int_{\partial\mathbb{B}^{n}}
\frac{-(2n-1)(1-|z|^{2})^{2n-2}[z_{k}|\zeta-z|^{2}+(1-|z|^{2})(z_{k}-\zeta_{k})]}
{|z-\zeta|^{4n}}f_{j}(\zeta) \,d\sigma(\zeta)
$$
which hold clearly for each $j, k\in\{1,\ldots,n\}$. Now, we introduce
$$\Gamma_{f_{j}}=\sum_{k=1}^{n}(f_{j}(z))_{z_{k}}\cdot
\theta_{k}+\sum_{k=1}^{n}(f_{j}(z))_{\overline{z}_{k}}\cdot\overline{\theta}_{k}.
$$
Then the classical Cauchy-Schwarz inequality yields

\vspace{5pt}

$\ds \frac{\big|\Gamma_{f_{j}}\big |^{2}}{(2n-1)^{2}(1-|z|^{2})^{4n-4}}$
\begin{eqnarray*}
&=& \left |\sum_{k=1}^{n}\int_{\partial\mathbb{B}^{n}}
\frac{[\overline{z}_{k}|\zeta-z|^{2}+(1-|z|^{2})(\overline{z}_{k}-\overline{\zeta}_{k})]\theta_{k}}
{|z-\zeta|^{4n}}f_{j}(\zeta)\,d\sigma(\zeta)\right .\\
&&~~~~~~~~~\left . +\sum_{k=1}^{n}\int_{\partial\mathbb{B}^{n}}
\frac{[z_{k}|\zeta-z|^{2}+(1-|z|^{2})(z_{k}-\zeta_{k})]\overline{\theta}_{k}}
{|z-\zeta|^{4n}}f_{j}(\zeta)\,d\sigma(\zeta)\right |^{2}\\
&\leq&4\left
[\int_{\partial\mathbb{B}^{n}}\frac{\big [|z||\zeta-z|^{2}+(1-|z|^{2})|\zeta-z|\big]|f_{j}(\zeta)|}
{|z-\zeta|^{4n}}\,d\sigma(\zeta)\right ]^{2}\\
&\leq&4\left  [\int_{\partial\mathbb{B}^{n}}\frac{\big [|z||\zeta-z|+(1-|z|^{2})\big ]^{2}}
{|z-\zeta|^{4n-2}}\,d\sigma(\zeta)\right ]\left
[\int_{\partial\mathbb{B}^{n}}\frac{|f_{j}(\zeta)|^{2}}{|z-\zeta|^{4n}}\,d\sigma(\zeta)\right
],
\end{eqnarray*}
whence

\vspace{5pt}

$\ds \frac{\big|\Lambda_{f}\big |^{2}}{(2n-1)^{2}(1-|z|^{2})^{4n-4}}$
\begin{eqnarray*}
&=&\frac{\max_{\theta\in\partial\mathbb{B}^{n}}\left(\sum_{j=1}^{n}
|\Gamma_{f_{j}}|^{2}\right)}{(2n-1)^{2}(1-|z|^{2})^{4n-4}}\\
&\leq& 4\left [\int_{\partial\mathbb{B}^{n}}\frac{[|z||\zeta-z|+(1-|z|^{2})]^{2}}
{|z-\zeta|^{4n-2}}\,d\sigma(\zeta)\right ]\left
[\int_{\partial\mathbb{B}^{n}}\frac{\sum_{j=1}^{n}|f_{j}(\zeta)|^{2}}
{|z-\zeta|^{4n}}\,d\sigma(\zeta)\right ]\\
&\leq&\frac{4M^{2}}{(1-|z|)^{2}(1-|z|^{2})^{2n-1}}\left
[\int_{\partial\mathbb{B}^{n}}\frac{(1+|z|)^{2}}
{|z-\zeta|^{4n-2}}\,d\sigma(\zeta)\right]\\
&\leq&\frac{4M^{2}(1+|z|)^{2}}{(1-|z|)^{2}(1-|z|^{2})^{2n-1}}\left
[\int_{\partial\mathbb{B}^{n}} \frac{1}{|z-\zeta|^{4n-2}}\,d\sigma(\zeta)\right]\\
&\leq& \frac{4M^{2}(1+|z|)^{2}}{(1-|z|)^{2}(1-|z|^{2})^{4n-2}}.
\end{eqnarray*}
Hence
$$|\Lambda_{f}|^2\leq \frac{4(2n-1)^2M^{2}}{(1-|z|)^{4}},
$$
from which the inequality \eqref{eq2} follows.
\epf

\begin{defn}
A matrix-valued function $A(z)=\big(a_{i,j}(z)\big)_{n\times n}$ is
called {\it $h$-harmonic} if each of its entries $a_{i,j}(z)$ is a
$h$-harmonic mapping from an open subset
$\Omega\subset\mathbb{C}^{n}$ into $\mathbb{C}$.
\end{defn}

As an application of Theorem \ref{thm1}, we get

\begin{lem}\label{lemr1}
Suppose that $A(z)=\big(a_{i,j}(z)\big)_{n\times n}$ is a matrix-valued
$h$-harmonic mapping of $\mathbb{B}^{n}(r)$ such that $A(0)=0$ and
$|A(z)|\leq M$ in $\mathbb{B}^{n}(r).$ Then
$$|A(z)|\leq M\left[1-\frac{(r-|z|)^{2n-1}}{(r+|z|)^{2n-1}}\right].
$$
\end{lem}
\bpf For an arbitrary $\theta=(\theta_1,\ldots,\theta_n)^{T}\in\partial\mathbb{B}^{n}$, we
introduce
$$ P_{\theta}(z)=A(z)\theta=(p_{1}(z),\ldots,p_{n}(z))
$$
and let $F_{\theta}(\zeta)=P_{\theta}(r\zeta)$ for $\zeta\in\mathbb{B}^{n}$. By  Theorem \ref{thm1},
we see that
$$\left|F_{\theta}(\zeta)-\frac{(1-|\zeta|)^{2n-1}}{(1+|\zeta|)^{2n-1}}F_{\theta}(0)\right|\leq
M\left[1-\frac{(1-|\zeta|)^{2n-1}}{(1+|\zeta|)^{2n-1}}\right], \quad \zeta\in\mathbb{B}^{n},
$$
which is equivalent to
$$|P_{\theta}(z)|\leq M\left[1-\frac{(r-|z|)^{2n-1}}{(r+|z|)^{2n-1}}\right], \quad z\in\mathbb{B}^{n}(r).
$$
The arbitrariness of $\theta$  yields the desired inequality.
\epf

We recall the following result which is crucial for the proof of our next theorem.

\begin{Lem}{\rm (\cite[Lemma 1]{HG1}~or \cite[Lemma 4]{LX})}\label{lem-Liu}
Let $A$ be an $n \times n$ complex $(real)$ matrix. Then for $\theta\in\partial \mathbb{B}^{n}$,
the inequality
%$$|A\theta|\geq\frac{|\det A|}{|A|^{n-1}}
%$$
$|A\theta|\geq|\det A|\, |A|^{1-n}$
holds.
\end{Lem}

\begin{thm}\label{thm3w}
Suppose that $f\in \mathcal{HB}_{n}(\alpha)$, $f(0)=0$, $\det
J_{f}(0)=1$ and $\|f\|_{\mathcal{HB}_{n}(\alpha)}\leq M,$ where
$M$ is a positive constant. Then $f$ is univalent in
$\mathbb{B}^{n}(\rho/2)$, where
\be\label{extra12}
\rho=\frac{3^{\alpha}}{(2M)^{2n}(3^{\alpha}+4^{\alpha})}.
\ee
Moreover, the range $f(\mathbb{B}^{n}(\rho/2))$ contains a
univalent ball $\mathbb{B}^{n}(R)$, where
$$R\geq\frac{\rho}{4M^{2n-1}}.
$$
\end{thm} \bpf
For $\zeta\in\mathbb{B}^{n}$, let
$F(\zeta)=2f(\frac{1}{2}\zeta)$. Then
$$| F_{\zeta}(\zeta)|+| F_{\overline{\zeta}}(\zeta)|\leq\frac{M}{\Big(1-\frac{|\zeta|^{2}}{4}\Big)^{\alpha}}
\leq\frac{4^{\alpha}}{3^{\alpha}}M
$$
which gives
$$ |F_{\zeta}(\zeta)-F_{\zeta}(0)|\leq |F_{\zeta}(\zeta)|+|F_{\zeta}(0)|
\leq\left(1+\frac{4^{\alpha}}{3^{\alpha}}\right)M.
$$
Lemma \ref{lemr1} implies that
\beq\label{eq-1-th3}
\nonumber |F_{\zeta}(\zeta)-F_{\zeta}(0)|&\leq&
\left(1+\frac{4^{\alpha}}{3^{\alpha}}\right)M\left[1-\frac{(1-|\zeta|)^{2n-1}}{(1+|\zeta|)^{2n-1}}\right]\\
\nonumber
&=&\frac{2M(3^{\alpha}+4^{\alpha})}{3^{\alpha}}\frac{\big(C_{2n-1}^{1}|\zeta|+C_{2n-1}^{3}|\zeta|^{3}+\cdots+
C_{2n-1}^{2n-1}|\zeta|^{2n-1}\big)}{(1+|\zeta|)^{2n-1}}\\ \nonumber
&\leq&
\frac{2^{2n-1}(3^{\alpha}+4^{\alpha})M}{3^{\alpha}(1+|\zeta|)^{2n-1}}|\zeta|\\
%\nonumber
&\leq&\frac{2^{2n-1}(3^{\alpha}+4^{\alpha})M}{3^{\alpha}}|\zeta|,
\eeq
where $C_{n}^{k}={n\choose k}$ ($k=1,2, \ldots, n$) denote the binomial coefficients. Similarly,
\be\label{eq-2-th3}
|F_{\overline{\zeta}}(\zeta)-F_{\overline{\zeta}}(0)|\leq
\frac{2^{2n-1}(3^{\alpha}+4^{\alpha})M}{3^{\alpha}}|\zeta|.
\ee
On the other hand, for
$\theta\in\partial\mathbb{B}^{n}$, we infer from (\ref{eqbe1}),
(\ref{eqbe2}) and  Lemma \Ref{lem-Liu} that
\be\label{eq-3-th3}
\lambda_{F}(0)\geq\frac{\det J_{F}(0)}{\Lambda_{F}^{2n-1}(0)}\geq\frac{1}{M^{2n-1}}.
\ee

In order to prove the univalence of $F$ in $\mathbb{B}^{n}(\rho)$,
we choose two distinct points $\zeta'$ and $\zeta''$ in
$\mathbb{B}^{n}(\rho)$ with $\zeta'-\zeta''=|\zeta'-\zeta''|\theta$,
and let $[\zeta',\zeta'']$ denote the line segment with endpoints $\zeta'$  and $\zeta''$,
where
$$\rho=\frac{3^{\alpha}}{(2M)^{2n}(3^{\alpha}+4^{\alpha})}.
$$
Set $d\zeta=( d\zeta_{1}, \ldots, d\zeta_{n})^T$ and $(d\overline{\zeta}=( d\overline{\zeta}_{1}, \ldots,
d\overline{\zeta}_{n})^T$.
Then we infer from  \eqref{eq-1-th3}, (\ref{eq-2-th3}) and (\ref{eq-3-th3}) that
\begin{eqnarray*}
|F(\zeta')-F(\zeta'')|&\geq&
\left|\int_{[\zeta',\zeta'']}F_{\zeta}(0)\,d\zeta+F_{\overline{\zeta}}(0)\,d\overline{\zeta}\right|\\
&& \hspace {1cm} - \left|\int_{[\zeta',\zeta'']}(F_{\zeta}(\zeta)-F_{\zeta}(0))\,d\zeta+(F_{\overline{\zeta}}(\zeta)-
F_{\overline{\zeta}}(0))\,d\overline{\zeta}\right|\\
&\geq&|F_{\zeta}(0)\theta+F_{\overline{\zeta}}(0)\overline{\theta}|\int_{[\zeta',\zeta'']}\,|d\zeta|-
\frac{2^{2n}(3^{\alpha}+4^{\alpha})M}{3^{\alpha}}\int_{[\zeta',\zeta'']}|\zeta|\,|d\zeta|\\
&>&|\zeta'-\zeta''|\left\{\lambda_{F}(0)
-\frac{2^{2n}(3^{\alpha}+4^{\alpha})M}{3^{\alpha}}\rho\right\}\\
&\geq&|\zeta'-\zeta''|\left\{\frac{1}{M^{2n-1}}
-\frac{2^{2n}(3^{\alpha}+4^{\alpha})M}{3^{\alpha}}\rho\right\}\\
&=&0,
\end{eqnarray*}
where $\theta=\frac{d\zeta}{|d\zeta|}$. Thus, $F$ is univalent in $\mathbb{B}^{n}(\rho)$
which is equivalent to saying that $f$ is univalent in $\mathbb{B}^{n}(\rho/2)$.

Furthermore, for each $z$ with $|\zeta|=\rho$, we have
\begin{eqnarray*}
|F(\zeta)-F(0)|&\geq&
\left|\int_{[0,\zeta]}F_{\zeta}(0)\,d\zeta+F_{\overline{\zeta}}(0)\,d\overline{\zeta}\right|\\
&& \hspace{1cm} -\left|\int_{[0,\zeta]}(F_{\zeta}(\zeta)-F_{\zeta}(0))\,d\zeta+(F_{\overline{\zeta}}(\zeta)-
F_{\overline{\zeta}}(0))\,d\overline{\zeta}\right|\\
&\geq&\rho\left\{\frac{1}{M^{2n-1}}-\frac{2^{2n-1}(3^{\alpha}+4^{\alpha})M\rho}{3^{\alpha}}\right\}\\
&=&\frac{\rho}{2M^{2n-1}} ~\mbox{ (by \eqref{extra12})}
\end{eqnarray*}
showing the range $f(\mathbb{B}^{n}(\rho/2))$ contains a univalent ball $\mathbb{B}^{n}(R)$, where
$R\geq\rho/(4M^{2n-1}). $ The proof of this theorem is complete.
\epf

\end{document}